\documentclass[12pt,reqno]{amsart}
\usepackage{enumerate, latexsym, amsmath, amsfonts, amssymb, amsthm, color}

\def\Z{\Bbb Z}

\def\Q{\Bbb Q}

\def\l{\left}
\def\r{\right}
\def\bg{\bigg}
\def\({\bg(}
\def\){\bg)}
\def\t{\text}
\def\f{\frac}
\def\ch{{\rm ch}}

\def\ls{\leqslant}
\def\gs{\geqslant}

\def\sm{\setminus}

\def\al{\alpha}

\def\ve{\varepsilon}

\theoremstyle{plain}
\newtheorem{theorem}{Theorem}

\newtheorem{conjecture}{Conjecture}
\theoremstyle{definition}

\theoremstyle{remark}
\newtheorem{remark}{Remark}

 \vspace{4mm}

\begin{document}

\hbox{Preprint}

\title
[{On sums and products in a field}]
{On sums and products in a field}

\author
[Guang-Liang Zhou and Zhi-Wei Sun] {Guang-Liang Zhou and Zhi-Wei Sun}

\address {(Guang-Liang Zhou) Department of Mathematics, Nanjing
University, Nanjing 210093, People's Republic of China}
\email{1064021981@qq.com}

\address{(Zhi-Wei Sun) Department of Mathematics, Nanjing
University, Nanjing 210093, People's Republic of China}
\email{zwsun@nju.edu.cn}

\keywords{Fields, rational functions, restricted sums, restricted products.
\newline \indent 2010 {\it Mathematics Subject Classification}. Primary 11D85; Secondary 11P99, 11T99.
\newline \indent The second author is the corresponding author. This research was supported by the Natural Science Foundation of China (grant 11571162) and the NSFC-RFBR Cooperation and Exchange Program (grant 11811530072).}

\begin{abstract} In this paper we study sums and products in a field. Let $F$ be a field with $\ch(F)\not=2$, where $\ch(F)$ is the characteristic of $F$.
For any integer $k\gs4$, we show that each $x\in F$ can be written as $a_1+\ldots+a_k$ with $a_1,\ldots,a_k\in F$ and $a_1\ldots a_k=1$ if $\ch(F)\not=3$, and that for any $\al\in F\sm\{0\}$ we can write each $x\in F$ as $a_1\ldots a_k$ with $a_1,\ldots,a_k\in F$ and $a_1+\ldots+a_k=\al$. We also prove that for any $x\in F$ and $k\in\{2,3,\ldots\}$ there are $a_1,\ldots,a_{2k}\in F$ such that $a_1+\ldots+a_{2k}=x=a_1\ldots a_{2k}$.
\end{abstract}
\maketitle

\section{Introduction}
\setcounter{lemma}{0}
\setcounter{theorem}{0}
\setcounter{corollary}{0}
\setcounter{remark}{0}
\setcounter{equation}{0}

Let $\Q$ be the field of rational numbers. In 1749 Euler showed that any $q\in\Q$ can be written as $abc(a+b+c)$ with $a,b,c\in\Q$;
equivalently, we can always write $x=-q\in\Q$ as $abcd$ with $a,b,c,d\in\Q$ and $a+b+c+d=0$.
Actually, Euler noted that the equation $abc(a+b+c)=q$ has the following rational parameter solutions:
\begin{align*}a=&\frac{6qst^3(qt^4-2s^4)^2}{(4qt^4+s^4)(2q^2t^8+10qs^4t^4-s^8)},
\\b=&\f{3s^5(4qt^4+s^4)^2}{2t(qt^4-2s^4)(2q^2t^8+10qs^4t^4-s^8)},
\\c=&\f{2(2q^2t^8+10qs^4t^4-s^8)}{3s^3t(4qt^4+s^4)}.
\end{align*}
The reader may consult N. D. Elkies's talk \cite{E} for a nice exposition of this curious discovery of Euler and its connection to modern topics like $K3$ surfaces. Elkies \cite{E} found that $abcd=x$ with $a+b+c+d=0$, where
\begin{gather*}a=\frac{(s^4+4x)^2}{2s^3(s^4-12x)},\ b=\f{2x(3s^4-4x)^2}{s^3(s^4+4x)(s^4-12x)},
\\c=\f{s(s^4-12x)}{2(3s^4-4x)},\ d=-\f{2s^5(s^4-12x)}{(s^4+4x)(3s^4-4x)}.
\end{gather*}

Let $F$ be a field. If $x=a_1\cdots a_k$ with $a_1,\ldots,a_k\in F$ and $a_1+\ldots+a_k=0$, then $a_1\ldots a_k$
is called a {\it balanced decomposition} of $x$ by A. A. Klyachko and A. N. Vassilyev \cite{KV}.
Unaware of Euler's above work in 1749, Klyachko and Vassilyev \cite{KV} showed that if $\ch(F)$ (the characteristic of $F$) is not two then for each $k=5,6,\ldots$
every $x\in F$ has a balanced decomposition $a_1\ldots a_k$ with $a_1,\ldots,a_k\in F$ and $a_1+\cdots+a_k=0$.
When $F$ is a finite field and $k>1$ is an integer, they determined completely when each $x\in F$ can be written as $a_1\ldots a_k$ with $a_1,\ldots,a_k\in F$ and $a_1+\cdots+a_k=0$. In 2016, Klyachko, A. M. Mazhuga and A. N. Ponfilenko \cite{KMP} proved that if $\ch(F)\not=2,3$ and $|F|\not=5$ then each $x\in F$ has a balanced decomposition $a_1a_2a_3a_4$ with $a_1,a_2,a_3,a_4\in F$ and $a_1+a_2+a_3+a_4=0$; in fact, for $x\in F\sm\{1/4,-1/8\}$ they found that $a(x)b(x)c(x)d(x)=x$ and $a(x)+b(x)+c(x)+d(x)=0$, where
\begin{gather*}a(x)=\f{2(1-4x)^2}{3(1+8x)},\ b(x)=-\f{1+8x}6,
\\c(x)=-\f{1+8x}{2(1-4x)},\ d(x)=\f{18x}{(1-4x)(1+8x)}.
\end{gather*}
This is much simpler than Euler's and Elkies' rational parameter solutions to the equation $abcd=x$ with the restriction $a+b+c+d=0$.

Motivated by the above work, we obtain the following new results.

\begin{theorem}\label{Th1.1} Let $F$ be any field with $\ch(F)\not=2$, and let $\alpha\in F\sm\{0\}$ and $k\in\{4,5,\ldots\}$. Then each $x\in F$
can be written as $a_1\ldots a_k$ with $a_1,\ldots,a_k\in F$ with $a_1+\ldots+a_k=\al$.
\end{theorem}
\begin{remark}\label{Rem1.1} This theorem with $\al=1$ implies that for any $x\in\Q$ and $k\in\{4,5,\ldots\}$ there are $a_1\ldots a_k\in\Q$ such that $a_1\ldots a_k(a_1+\ldots+a_k)=x$, this extension of Euler's work was recently asked by D. van der Zypen \cite{Z}.
\end{remark}

\begin{theorem}\label{Th1.2} Let $F$ be a filed with $\ch(F)\not=2$, and let $k\gs4$ be an integer.

{\rm (i)} Let $x\in F$. If $\ch(F)\not=3$ or $x\not=0$, then we can write $x$ as $a_1+\ldots+a_k$ with $a_1,\ldots,a_k\in F$ and $a_1\ldots a_k=-1$.

{\rm (ii)} If $\ch(F)\not=3$, then any $x\in F$ can be written as $a_1+\ldots+a_k$ with $a_1,\ldots,a_k\in F$ and $a_1\ldots a_k=1$.
\end{theorem}
\begin{remark}\label{Rem1.2} It seems that there are no $a,b,c\in\Q$ with $a+b+c=1=abc$.
We guess that the condition $\ch(F)\not=3$ in part (iii) of Theorem 1.2 can be removed.
For the finite field $\mathbb F_3=\{0,\pm1\}$ of order $3$, clearly
\begin{align*}0=&1+1-1-1\quad\t{with}\ 1\times1\times(-1)\times(-1)=1,
\\1=&1+1+1+1\quad\t{with}\  1\times1\times1\times1=1,
\\-1=&-1-1-1-1\quad\t{with}\  (-1)\times(-1)\times(-1)\times(-1)=1.
\end{align*}
\end{remark}

Let $F$ be any field and $k$ be a positive integer. Clearly, any $x\in F$ can be written as $a_1+\ldots +a_{2k+1}$
with $a_1,\ldots,a_{2k+1}\in F$ and $a_1\ldots a_{2k+1}=(-1)^kx$; in fact, $x+k(1-1)=x$ and $x\times1^k\times(-1)^k=(-1)^kx$.
If $a^2=-1$ for some $a\in F$, then any $x\in F$ can be written as $a_1+\ldots +a_{2k+1}$
with $a_1,\ldots,a_{2k+1}\in F$ and $a_1\ldots a_{2k+1}=(-1)^{k-1}x$; in fact, $x+(a-a)+(k-1)(1-1)=x$ and
$$x\times a\times(-a)\times 1^{k-1}\times(-1)^{k-1}=(-1)^{k-1}x.$$

\begin{theorem}\label{Th1.3} Let $F$ be a field with $\ch(F)\not=2$, and let $k\gs2$ be an integer.
Then, for any $x\in F$ there are $a_1,\ldots,a_{2k}\in F$ such that $a_1+\ldots+a_{2k}=x=a_1\ldots a_{2k}$.
\end{theorem}
\begin{remark}\label{Rem1.3} If $F$ is a field with $\ch(F)\not=2$, and $k\gs2$ is an integer, then for any $x\in F$,
by Theorem 1.3 there are $a_1,\ldots,a_{2k}\in F$ with $a_1+\ldots+a_{2k}=-x=a_1\ldots a_{2k}$, hence
$$(-a_1)+\ldots+(-a_{2k})=x\quad\t{and}\quad (-a_1)\ldots(-a_{2k})=-x.$$
\end{remark}

Motivated by Theorem 1.3 and Remark 1.3, we propose the following conjecture based on our computation.

\begin{conjecture}\label{Conj1.1} Let $F$ be any field with $\ch(F)\not=2,3$. Then, for any $x\in F$ there are $a,b,c,d\in F$
such that $a+b+c+d-1=x=abcd$.
\end{conjecture}

For example, in any field $F$ with $\ch(F)\not=2,3$, we have
$$-2+\f92-\f23+\f16-1=1=(-2)\times\f 92\times\l(-\f23\r)\times\f16.$$

\section{Proofs of Theorems 1.1-1.3}
\setcounter{lemma}{0}
\setcounter{theorem}{0}
\setcounter{corollary}{0}
\setcounter{remark}{0}
\setcounter{equation}{0}

\medskip
\noindent{\it Proof of Theorem} 1.1. We distinguish three cases.

{\it Case} 1. $k=4$.

If each $q\in F$ can be written as $abcd$ with $a,b,c,d\in F$ and $a+b+c+d=1$, then
for any $x\in F$ we can write $x/\al^4=abcd$ with $a,b,c,d\in F$ and $a+b+c+d=1$ and hence $x=(a\al)(b\al)(c\al)(d\al)$
with $a\al+b\al+c\al+d\al=\al$. So, it suffices to work with $\al=1$.

Let $x\in F$ with $x\not=\pm1$. Define
$$a(x)=-\frac{(1-x)^2}{2(1+x)},\ b(x)=\f{1+x}2,\ c(x)=\f{1+x}{1-x},\ d(x)=\f{4x}{x^2-1}.$$
It is easy to verify that
$$a(x)b(x)c(x)d(x)=x\ \ \text{and}\ \ a(x)+b(x)+c(x)+d(x)=1.$$

For $x=-1$, we note that
$$-1=\f12\times\f12\times2\times(-2)\quad\t{with}\ \f12+\f12+2-2=1.$$
For $x=1$, if $\ch(F)=3$ then
$$1=1\times1\times1\times1\quad\t{with}\ 1+1+1+1=1,$$
if $\ch(F)\not=3$ then
$$1=\f32\times\l(-\f32\r)\times\l(-\f13\r)\times\f43
\ \ \ \text{with}\ \f32-\f32-\f13+\f43=1.$$
This proved Theorem 1.1 for $k=4$.

{\it Case} 2. $k=5$.

As $\ch(F)\not=2$, we have $\al-\ve\not=0$
for some $\ve\in\{\pm1\}$. Let $x\in F$. By Theorem 1.1 for $k=4$, we can write $\ve x$ as $abcd$ with $a,b,c,d\in F$ and $a+b+c+d=\al-\ve$.
Hence $x=abcd\ve$ with $a+b+c+d+\ve=\al$.
So Theorem 1.1 also holds for $k=5$.

{\it Case} 3. $k\gs6$.

Let $x\in F$. If $k$ is even, then by Theorem 1.1 for $k=4$ there are $a,b,c,d\in F$ with $a+b+c+d=\al$
such that $abcd=(-1)^{(k-4)/2}x$, hence
$$x=abcd\times1^{(k-4)/2}\times(-1)^{(k-4)/2}$$
with $$a+b+c+d+\f{k-4}2(1-1)=\al.$$
When $k$ is odd,
by Theorem 1.1 for $k=5$ there are $a,b,c,d,e\in F$ with $a+b+c+d+e=\al$
such that $abcde=(-1)^{(k-5)/2}x$, hence
$$x=abcde\times1^{(k-5)/2}\times(-1)^{(k-5)/2}$$
with $$a+b+c+d+e+\f{k-5}2(1-1)=\al.$$

Combining the above, we have completed the proof of Theorem 1.1. \qed

\medskip
\noindent{\it Proof of Theorem} 1.2. For any $m\in\Z$, if each $x\in F$ can be written as $a+b+c+d$ with $a,b,c,d\in F$ and $abcd=m$,
then for any $x\in F$ and $k\in\{4,5,\ldots\}$ there are $a_1,a_2,a_3,a_4\in F$ such that $a_1+a_2+a_3+a_4=x-(k-4)$ and $a_1a_2a_3a_4=m$, hence $a_1+\ldots+a_k=x$ and $a_1\ldots a_k=m$, where $a_j=1$ for $4<j\ls k$. Thus it suffices to show parts (i) and (ii) in the case $k=4$.

(i) For $x\in F\sm\{-1,-3\}$, we define
$$a(x)=\f{(x+1)^2}{2(x+3)},\ b(x)=\f{x+3}2,\ c(x)=-\f{x+3}{x+1},\ d(x)=\f4{(x+1)(x+3)},$$
and it is easy to verify that
$$a(x)b(x)c(x)d(x)=-1\quad\t{and}\quad a(x)+b(x)+c(x)+d(x)=x.$$
Observe that
$$-1=2-2-\f12-\f12\quad\t{with}\ 2\times(-2)\times\l(-\f12\r)\times\l(-\f12\r)=-1.$$
If $\ch(F)\not=3$, then
$$-3=\f23-\f23-\f32-\f32\quad\t{with}\ \f23\times\l(-\f23\r)\times\l(-\f32\r)\times\l(-\f32\r)=-1.$$
When $\ch(F)=3$, the element $-3$ of $F$ is just zero.
This concludes the proof of Theorem 1.2(i).

(ii) Suppose that $\ch(F)\not=3$.  Clearly, $0=1+1-1-1$ with $1\times1\times(-1)\times(-1)=1$.
If $\ch(F)=5$ and $x\in F\sm\{0\}$, then
$$x=x^{-1}-x^{-1}-2x-2x\quad\t{and}\quad x^{-1}(-x^{-1})(-2x)(-2x)=-4=1.$$

Now we consider the case $\ch(F)\not=5$. For $x\in F\sm\{-3/2,-9\}$, we define
\begin{gather*}a(x)=-\f{2(x+9)^2}{15(2x+3)},\ \ b(x)=\f{8(2x+3)}{15},
\\ c(x)=-\f{2x+3}{4(x+9)},\ \ d(x)=\f{25}{4(x+9)(2x+3)},
\end{gather*}
and it is easy to verify that
$$a(x)b(x)c(x)d(x)=1\quad\t{and}\quad a(x)+b(x)+c(x)+d(x)=x.$$
Note that
$$-\f32=1-1-2+\f12\quad\t{with}\ 1\times(-1)\times(-2)\times\f12=1$$
and
$$-9=-4-\f92-\f13-\f16\quad\t{with}\ (-4)\times\l(-\f92\r)\times\l(-\f13\r)\times\l(-\f16\r)=1.$$
So Theorem 1.2(ii) also holds.

In view of the above, the proof of Theorem 1.2 is now complete. \qed

\medskip
\noindent{\it Proof of Theorem} 1.3.
We first handle the case $k=2$. For $x\in F\sm\{\pm1\}$, we define
$$a(x)=\f{(x+1)^2}{2(x-1)},\ b(x)=\f{x-1}2,\ c(x)=\f{1-x}{1+x},\ d(x)=\f{4x}{1-x^2},$$
and it is easy to verify that
$$a(x)b(x)c(x)d(x)=x=a(x)+b(x)+c(x)+d(x).$$
Clearly,
$$-1=-\f12-\f12+2-2\ \ \ \t{with}\ -\f12\times\l(-\f12\r)\times2\times(-2)=-1.$$
If $\ch(F)\not=3$, then
$$1=\f32-\f32-\f13+\f43\ \ \ \t{with}\ \f32\times\l(-\f32\r)\times\l(-\f13\r)\times\f43=1.$$
When $\ch(F)=3$, we have
$$1=1+1+1+1\quad\t{with}\ 1\times1\times1\times1=1.$$
This proves Theorem 1.3 for $k=2$.

Now we consider the case $k\gs3$. By Theorem 1.3 for $k=2$, there are $a,b,c,d\in F$ such that $a+b+c+d=(-1)^kx=abcd$.
Thus $$(-1)^ka+(-1)^kb+(-1)^kc+(-1)^kd+(k-2)(1-1)=x$$ and
$$(-1)^ka\times(-1)^kb\times(-1)^kc\times(-1)^kd\times1^{k-2}\times(-1)^{k-2}=abcd(-1)^k=x.$$
This proves Theorem 1.3 for $k\gs3$.

By the above, we have completed the proof of Theorem 1.3. \qed

\setcounter{conjecture}{0} \end{document}